\documentclass[12pt]{amsart}
\usepackage{amsfonts}
\usepackage{amsmath}
\usepackage{amsxtra}
\usepackage{amssymb,latexsym}
\usepackage[mathcal]{eucal}
\usepackage{amscd}

\newtheorem{theo}{{\bfseries Theorem}}[section]
\newtheorem{prop}[theo]{{\bfseries Proposition}}

\newtheorem{cor}[theo]{{\bfseries Corollary}}
\newtheorem{df}[theo]{{\bfseries Definition}}

\newtheorem{ex}[theo]{{\bfseries Example}}

\def \R {\mathbb R}

\def \O {\mathcal O}

\def \S {\mathcal S}

\def \d {\delta}

\def \l {\lambda}

\usepackage{amssymb,latexsym}
\usepackage[mathcal]{eucal}
\usepackage{amscd}
\usepackage{amsmath}
\usepackage{graphicx}
\usepackage{amsfonts}
\usepackage{amscd}

\numberwithin{equation}{section}

\begin{document}

\title[Smale Strategies for the N-Person IPD]{Smale Strategies for the \\ N-Person Iterated Prisoner's Dilemma }
\vspace{1cm}
\author{Ethan Akin}
\address{ Mathematics Department \\
    The City College \\
    137 Street and Convent Avenue \\
    New York City, NY 10031, USA }
\email{ethanakin@earthlink.net}

\author{Slawomir Plaskacz}
\address{Faculty of Mathematics and Computer Science\\
   N. Copernicus University in Torun }
\email{plaskacz@mat.umk.pl}

\author{Joanna Zwierzchowska}
\address{Faculty of Mathematics and Computer Science\\
   N. Copernicus University in Torun }
\email{joanna.zwierzchowska@mat.umk.pl}

    \vspace{.5cm}
\date{February, 2018}

\begin{abstract}
Adapting methods introduced by Steven Smale, we describe good strategies
for a symmetric version of the Iterated Prisoner's Dilemma with
$n$ players.
\end{abstract}



 \maketitle

 \section{Introduction}
 \vspace{.5cm}

 In \cite{Sm} Smale introduced an approach to strategy for the Iterated Prisoner's Dilemma which was different from the popular Markov chain
 methods. He suggested using as data the current time average payoff to the players rather than using the results of the most recent round of play.
 Smale's results were extended in \cite{A17}. Here we apply these methods to a symmetric $n$ player version of the Prisoner's Dilemma. The goal is
 to describe good plans which stabilize the cooperative outcome where each player receives (in long-term average) the payoff $p_n$ obtained
 when all $n$ players cooperate. A strategy is a choice of initial play together with a plan responding to previous play.

 We describe good plans with the following properties.
 \begin{itemize}
 \item If all $n$ players eventually use good plans, then from any initial position, convergence to  cooperation is achieved.

 \item Suppose that for some $k$ with $1 \leq k < n$,  players $j$ with $j = 1,\dots,k$ eventually use good plans and suppose that
 $X = (x_1, \dots, x_k, z_{k+1}, \dots, z_n)$ is a limit point for the sequence of time averages of the payoffs. Let $\bar x$ be the mean of
 the good player payoffs and $\bar z$ the mean of the remaining, dissenting, player payoffs.  The only way it can occur that $\bar z \geq p_n$
 is if $X = (p_n, \dots, p_n)$, i.e. the cooperative payoff.  Otherwise, while $\bar z > \bar x$, i.e. on average the dissenters do better than
 the good players, their  average payoff remains below $p_n$.

 If $k < n-1$ it may happen that, while the mean $\bar z < p_n$, some of the dissenters do get a payoff larger than $p_n$.  However,
 if $z_{k+1} \geq p_n$, then if $\tilde z$ is the mean of the remaining dissenters, i.e. the average of $z_{k+2}, \dots,z_n$ then
 $\tilde z < \bar x$. That is, the remaining players on average do worse than the good players and thus it pays for some of them to
 switch to good plans.
 \end{itemize}

Smale strategies for the Iterated Prisoner's Dilemma with $n$ players were described by Behrstock, Benaim and Hirsch in \cite{BBH}.
There is overlap between our results and theirs, but because they were considering a variety of cases, their exposition is
 more complicated.  Markov strategies for the Iterated Prisoner's Dilemma with $n$ players were described in \cite{HWTN14} and \cite{HWTN15}.

 \section{The Game}
  \vspace{.5cm}

 We consider an n-player symmetric game.  Each player has a choice of two strategies $c$ or $d$ (\emph{cooperate} or \emph{defect}).
 When no players cooperate, each obtains a payoff $r_0$ and when all n players cooperate, each obtains $p_n$. For $k = 1, \dots, n -1$
 when $k$ players cooperate, each cooperator receives $p_k$ and each defector receives $r_k$.

 Our first assumption is monotonicity of the payoffs.
 \begin{equation}\label{1.01}
 r_0 < r_1 < \dots < r_{n-1}, \quad \text{and} \quad   p_1 < p_2 < \dots < p_{n}.
 \end{equation}
 That is, the payoffs to both defectors and cooperators increase as the number of cooperators increase.

 Our second assumption concerns the advantage of defection.
 \begin{equation}\label{1.02}
 \text{For \ \ } k = 1, \dots, n, \quad p_k < r_{k-1}.
 \end{equation}
 This says that any cooperator would do better by switching to defection, provided the choices of the other players remain fixed.
 Thus, the strategy of defection \emph{dominates} cooperation regardless of the choices of the other players.

 For our third assumption, we define the population mean payoff, or just \emph{mean payoff}, $m_k$ by
 \begin{align}\label{1.03}
 \begin{split}
 m_0 \ = \ r_0, \quad m_n \ = \ p_n,& \\
 m_k \ = \ \frac{1}{n}[ k \cdot p_k \ + \ (n-k) \cdot r_k]& \quad \text{for \ \ }  k = 1, \dots, n-1.
 \end{split}
 \end{align}

 From (\ref{1.01}) and (\ref{1.02}) it is clear than
 \begin{equation}\label{1.03a}
 p_k < m_k < r_k \quad \text{for \ \ }  k = 1, \dots, n-1.
 \end{equation}

 Our third assumption is monotonicity of the mean payoffs.
  \begin{equation}\label{1.04}
 m_0 < m_1 < \dots < m_{n-1}  < m_{n}. \hspace{2cm}
 \end{equation}
 If there are currently $k-1$ cooperators and one defector switches to cooperation, then he bears a cost of $r_{k-1} - p_k$.
 This assumption says that the aggregate increase to the other players exceeds this cost, i.e.
  \begin{align}\label{1.05}
 \begin{split}
\ r_0 - p_1 < (n-1)\cdot(r_1 - r_0), \quad &r_{n-1} - p_n < (n-1)\cdot(p_n - p_{n-1}),\\
r_{k-1} - p_k <%PZ\ = \
(k-1) \cdot [(p_k  -  &p_{k-1}) \ + \ (n-k) \cdot (r_k - r_{k-1})], \\
 &\text{for} \quad k = 2, \dots, n-1.
 \end{split}
 \end{align}

 In the case when $n = 2$ the game is the classic Prisoner's Dilemma. In the notation of \cite{A17}
\begin{equation}\label{1.06}
 r_0 = P, \ r_1 = T, \ p_1 = S, \ p_2 = R,
 \end{equation}
 with the assumptions
 \begin{equation}\label{1.07}
 S < P < \frac{1}{2}(T + S) < R < T.
 \end{equation}

%EJAnew
 \begin{ex}\label{exsimple}  The values $(p_1,r_0,p_2,r_1, \dots) = (0,1,2,3, \dots)$ for $n \geq 2$ is an example. \end{ex}

 {\bfseries Proof:}  $p_k = 2k - 2, r_k = 2k + 1$ and so
 $$n \cdot m_k \ = \ k p_k + (n - k) r_k \ = \ k (2k - 2) + (n - k)(2k + 1) \ = \ (2n - 3) k + n.$$
  Thus,
  \begin{equation}\label{1.05a}
  m_k \ = \ (2 - \frac{3}{n}) k + 1.  \hspace{3cm}
  \end{equation}
  which is positive and increasing in $k$ for $n \geq 2$.

  $\Box$ \vspace{.5cm}
%
% As an example, $r_k = k + 1.5, $ for $k = 0, \dots n-1$ and $p_k = k$ for $k = 1, \dots, n$ satisfy the three assumptions.

 \vspace{1cm}

  \section{The Outcome and Strategy Regions}
   \vspace{.5cm}

   We consider repeated play.  After the $T^{th}$ round the payoff to the n players is a vector $S^T$ in $\R^n$. There are $2^n$
   possible payoff vectors according to the choices made by the players. For example, if in round $T$ a single player cooperated,
   then the vector $S^T$ has a single entry $p_1$ and n-1 entries $r_1$. There are n such payoff vectors according to which player
   cooperated.

   For a single payoff after $T$ rounds we use the time-average of the payoff vectors:
   \begin{equation}\label{2.01}
   s^T \ = \ \frac{1}{T} \sum_{t=1}^T \ S^t.
   \end{equation}

  Observe that
\begin{equation}\label{2.02}
s^{T+1} \quad = \quad \frac{T}{T+1} s^T \ + \ \frac{1}{T+1} S^{T+1}.
\end{equation}
and so
\begin{equation}\label{2.03}
s^{T+1} - s^T \quad = \quad  \frac{1}{T+1} (S^{T+1} - s^T).
\end{equation}

The sequence $\{ s^T \}$ lies in the \emph{outcome set} $\O \subset \R^n$ which is a convex hull of the
$2^n$ payoff vectors. The sequence need not converge in
$\O$ but the set $\Omega$ of limit points is always a nonempty, compact, connected subset of $\O$, see, e.g. \cite{A17} Proposition 2.1.

The approach introduced by Smale in \cite{Sm} is for each player to use the current time average $s^T$ as data  in order to determine
his choice in round $T+1$. In this case, the symmetry of the game suggests that it would be sufficient for each player to keep track of
how he is doing and to compare it merely with the time average of the mean payoff. That is, define $\bar \pi : \R^n \to \R$ by
$\bar \pi(x_1,\dots,x_n) = \frac{1}{n} \sum_{j = 1}^n \ x_j$. For player $j$ the linear map $\Pi_j = \pi_j \times \bar \pi : \R^n \to \R^2$ maps
the outcome set $\O$ onto the \emph{strategy set} $\S \subset \R^2$ which is the convex hull of the set of $2n$ points
\begin{equation} \label{2.04}
\{ (r_0,r_0), (p_n,p_n) \} \cup \{ (p_k,m_k), (r_k,m_k) : k = 1, \dots, n-1 \}
\end{equation}
Thus, $\Pi_j$ maps $s^T$ to a pair $(x,y) \in \R^2$ with $x$ the $j$ player's current average payoff and with $y$ the
current average population mean.

Notice that from  (\ref{1.04}) %PZ(\ref{1.05})
and convexity it follows that $(r_0,r_0)$ is the unique point of
$\S$ with $y$ coordinate less than or equal to $r_0$ and
$(p_n,p_n)$  is the unique point of $\S$ with $y$ coordinate greater  than or equal to $p_n$. Hence, these two points are extreme points of $\S$.

We will call $\{ (p_k,m_k) : k = 1, \dots, n \}$ the \emph{cooperation points} and
$\{ (r_k,m_k) : k = 0, \dots, n-1 \}$ the \emph{defection points} of $\S$.

\begin{df}\label{def2.1} For a line $\ell$ in $\R^2$ with equation $y = \l x + b$ we call the map $L : \R^2 \to \R$ given by
$(x,y) \mapsto y - \l x - b$ the \emph{associated affine map}. We call $\ell$ a \emph{separation line} for the game if
each cooperation point is on or above $\ell$ and each defection point is on or below $\ell$, i.e. $L$ is non-negative on the
cooperation points and non-positive on the defection points.\end{df}
   \vspace{.5cm}

   \begin{prop}\label{prop2.2} If a line $\ell$ with equation $y = \l x + b$ is a separation line, then $0 \leq \l \leq 1$ and
   if $n > 2$ then $0 < \l$. \end{prop}

{\bf Proof:} By the Intermediate Value Theorem, a separation line intersects the intervals $[(p_1,m_1),(r_0,r_0)]$
and $[(p_n,p_n),(r_{n-1},m_{n-1})]$. Hence, the steepest possible separation line connects $(r_0,r_0)$ and $(p_n,p_n)$ with slope 1.
On the other hand, the slope must be at most that of the line connecting $(p_1,m_1)$ with $ (r_{n-1},m_{n-1})$. This line has
positive slope if $n > 2$ and slope 0 if $n = 2$.

$\Box$ \vspace{.5cm}

%EJAnew

We will call the line through $(r_0,r_0)$ and $(p_n,p_n)$ the \emph{diagonal line} or just the diagonal.  It is the line with equation $y = x$ and
is the unique separation line with slope equal to $1$.

All our results are based on the following.

\begin{prop}\label{prop2.3} Assume that $X$ is a limit point of the sequence $\{ s^T \}$, i.e. $X \in \Omega \subset \O$. Let $\ell$ be a line
in $\R^2$ with finite slope.

(a) Assume that every defection point lies on or below $\ell$. If from some time onward, player $j$ defects whenever
$\Pi_j(s^T)$ lies above $\ell$, then $\Pi_j(X)$ lies on or below $\ell$.

(b) Assume that every cooperation point lies on or above $\ell$. If from some time onward, player $j$ cooperates whenever
$\Pi_j(s^T)$ lies below $\ell$, then $\Pi_j(X)$ lies on or above $\ell$.
\end{prop}

{\bf Proof:} Assume that $y = \l x + b$ is the equation of $\ell$ and that $L$ is the associated affine map. Let $M$ be a positive
number greater than the maximum of $L$ on $\S$.

(a) Assume that after some time $T_0 \geq 1$ player $j$ defects when $L(\Pi_j(s^T)) > 0$. Let $M_0 = T_0 \cdot M$. We show that for all $T \geq T_0$,
$L(\Pi_j(s^T)) \leq M_0/T$. We proceed by induction. Observe that the result holds for $T = T_0$.

If $L(\Pi_j(s^T)) > 0$ then player $j$ defects on the next round and so $L(\Pi_j(S^{T+1})) \leq 0$. Hence, because $L$ is an affine map,
(\ref{2.02}) implies that
\begin{equation}\label{2.05}
\begin{split}
L(\Pi_j(s^{T+1})) = \frac{T}{T+1} L(\Pi_j(s^T)) \ + \ \frac{1}{T+1} L(\Pi_j(S^{T+1})) \\
\leq \frac{T}{T+1} L(\Pi_j(s^T)) \leq \frac{T}{T+1} M_0/T = M_0/(T + 1),
\end{split}
\end{equation}
by induction hypothesis.

On the other hand, if $L(\Pi_j(s^T)) \leq 0$
\begin{equation}\label{2.06}
\begin{split}
L(\Pi_j(s^{T+1})) = \frac{T}{T+1} L(\Pi_j(s^T)) \ + \ \frac{1}{T+1} L(\Pi_j(S^{T+1})) \\
\leq \frac{1}{T+1} L(\Pi_j(S^{T+1})) \leq \frac{1}{T+1} M  \leq M_0/(T + 1).
\end{split}
\end{equation}

It follows that $\lim \sup \{ L(\Pi_j(s^T)) \} \leq 0$ and so any limit point lies on or below $\ell$.

(b) The proof is completely analogous and left to the reader. It is convenient to use $-L$ instead of $L$.

$\Box$ \vspace{1cm}

  \section{Simple Smale Plan}
   \vspace{.5cm}

\begin{df}\label{def3.1} Let $\ell$ be a separation line for the game %PZ we added -for the game-
with positive slope. The \emph{simple Smale plan} for player $j$
defects on round $T + 1$ when $\Pi_j(s^T)$ is above $\ell$ and cooperates when $\Pi_j(s^T)$ is on or below $\ell$.
\end{df} \vspace{.5cm}

We use the term ``plan'' because a strategy includes a choice of initial play as well as a plan responding to previous outcomes.
Observe that the demand that the slope be positive is automatically satisfied if $n > 2$ by Proposition \ref{prop2.2}.

%EJA
We will say that a player \emph{eventually uses a simple Smale plan} when there is a fixed simple Smale plan which is used by the player from
some time onward.

\begin{theo}\label{theo3.2} Let $\Omega$ be the set of limit points of the sequence $\{ s^T \}$.
If player $j$ eventually uses the simple Smale plan with separation line $\ell$, then  $\Pi_j(\Omega) \subset \ell$.\end{theo}

{\bf Proof:} This immediately follows from Proposition \ref{prop2.3}.

$\Box$ \vspace{.5cm}

\begin{cor}\label{cor3.3} For $j = 1, \dots, n$ let $\ell_j$ be a separation line with positive slope and assume that at least one
has slope less than 1. There exists a unique point $X \in \O$ such that if each player $j$ eventually uses  the simple Smale plan
associated with $\ell_j$, then the sequence $\{ s^T \}$ converges in $\O$ to $X$.

If each $\ell_j$ passes through $(p_n,p_n)$, then $X = (p_n, \dots, p_n)$ and so we obtain convergence to the position of complete cooperation.

\end{cor}

{\bf Proof:} Since  line $\ell_j$ has positive slope we can write its equation as $y = \l_j(x - a_j)$.  Now let $(x_1, \dots, x_n)$
be an arbitrary limit point of the sequence $\{ s^T \}$. If player $j$ uses the $\ell_j$ strategy, then with $\bar x = \bar \pi(x_1,\dots,x_n)$
Theorem \ref{theo3.2} implies that $(x_j, \bar x) = \Pi_j(x_1,\dots,x_n)$ lies on $\ell_j$ and so $\bar x = \l_j(x_j - a_j)$.

If this holds for all $j$ then $(\l_j)^{-1} \bar x = x_j - a_j$. Let $\bar a = \bar \pi(a_1,\dots,a_n)$ and
let $\hat \l$ be the harmonic mean of the $\l_j$'s, i.e. \\ $\hat \l = [\bar \pi((\l_1)^{-1},\dots,(\l_n)^{-1})]^{-1}$.
It is clear that $ \bar x = \hat \l (\bar x - \bar a)$. Since not all the $\l_j = 1$, Proposition \ref{prop2.2} implies that
$0 < \hat \l < 1$ and so the line $y = \hat \l (x - \bar a)$ intersects the line $y = x$ at a unique point $(\bar x, \bar x)$.
That is, the mean, $\bar x$, is determined by the lines $\ell_j$. Since $\bar x = \l_j(x_j - a_j)$, all of the $x_j$'s are determined as
well.

Thus, the limit point $X = (x_1,\dots,x_n)$ is uniquely determined by the equations $ \bar x = \hat \l (\bar x - \bar a)$
and $\bar x = \l_j(x_j - a_j)$ for $j = 1, \dots, n$.

Since the sequence $\{ s^T \}$ has a unique limit point, compactness implies convergence.

Now assume that $(p_n,p_n) \in \ell_j$ for all $j$.   Hence, $(\l_j)^{-1} p_n = p_n - a_j$ for all $j$ and so, averaging, we have
$ p_n = \hat \l (p_n - \bar a)$. Uniqueness in the above argument implies $\bar x = p_n$. The point $X = (p_n, \dots, p_n)$ is the
unique point of $\O$ with $\bar \pi(X) = p_n$.

$\Box$ \vspace{.5cm}

%EJA
{\bfseries Remark:} We can characterize the limit point $(x_1,\dots,x_n)$ as follows. It is the unique point such that there exists $y \in \R$
so that $X_k = (x_k,y) \in \ell_k$ and $\bar X = \frac{1}{n} \sum_{k=1}^n X_k$ is on the diagonal line with $y = x$.
\vspace{.5cm}

\begin{df}\label{def3.4} Let $\ell$ be a separation line with positive slope. We call simple Smale plan a
\emph{good simple Smale plan}
when it is associated with a separation line $\ell$ which passes through $(p_n,p_n)$ and which has slope $\l$
satisfying $\frac{n-1}{n} < \l < 1$.
\end{df} \vspace{.5cm}

That is, the equation of a separation line for a good simple Smale plan can be written $p_n - y = \l(p_n - x)$ with $\frac{n-1}{n} < \l < 1$.

It follows from Corollary \ref{cor3.3} that if each player eventually uses a good simple Smale plan then we achieve convergence to complete
cooperation. Smale calls this \emph{stability}.  However, we are also interested in what happens when only some of the players are good.

Suppose that for $j = 1,\dots, k$ player $j$ uses a good simple Smale plan associated with separation line $\ell_j$ having slope $\l_j$.
Suppose further that $X = (x_1,\dots,x_k,z_{k+1},\dots,z_n)$ is a limit point for $\{ s^T \}$. We assume that $1 \leq k \leq n-1$. We will call
the players $j$ with $j = k+1, \dots, n$ the \emph{dissenting players}.

Define
\begin{equation}\label{3.01}
\begin{split}
\bar x = \frac{1}{k} \sum_{j = 1}^k x_j, \quad \bar z = \frac{1}{n - k} \sum_{j = k+1}^n z_j, \hspace{2cm}\\
y = \frac{1}{n} [\sum_{j = 1}^k x_j +  \sum_{j = k+1}^n z_j] = \frac{1}{n} [k \bar x + (n-k) \bar z], \\
\hat \l = [\frac{1}{k} \sum_{j = 1}^k (\l%PZ\ell
_j)^{-1})]^{-1}.  \hspace{3cm}
\end{split}
\end{equation}

Thus, $\bar x$ is the mean payoff among the good players, $\bar z$ is the mean payoff among the dissenters and $y$ is the mean for the
entire population.

From Theorem \ref{theo3.2} it follows that $(x_j,y) = \Pi_j(X) \in \ell_j$ and so $p_n - y = \l_j(p_n - x_j)$ for $j = 1, \dots, k$.
Dividing by $\l_j$ and averaging as before we obtain $p_n - y = \hat \l(p_n - \bar x)$. Hence,
\begin{equation}\label{3.02}
\frac{1}{n} [k (p_n - \bar x) + (n-k) (p_n - \bar z)] = \hat \l(p_n - \bar x).
\end{equation}
We can rewrite this as
\begin{equation}\label{3.03}
\begin{split}
(n-k) (p_n - \bar z) = (n \hat \l - k) (p_n - \bar x), \qquad \text{and} \\
\bar z = \bar x + \frac{n(1 - \hat \l)}{n - k}(p_n - \bar x) \hspace{2cm}
\end{split}
\end{equation}

Now assume that $k \leq n - 2$ and suppose that $z_{k+1} \geq p_n$ and so $p_n - z_{k+1} \leq 0$. Define
\begin{equation}\label{3.04}
\tilde z = \frac{1}{n - k - 1} \sum_{j = k+2}^n z_j.
\end{equation}
 So we have
 \begin{equation}\label{3.05}
\frac{1}{n} [k (p_n - \bar x) + (n-k-1) (p_n - \tilde z)] \geq \hat \l(p_n - \bar x).
\end{equation}
We can rewrite this as
 \begin{equation}\label{3.06}
 \tilde z \leq %PZ=
 \bar x - \frac{1 - n(1 - \hat \l)}{n - k - 1}(p_n - \bar x).
 \end{equation}

 Now we apply all this.

 \begin{theo}\label{theo3.5} With $1 \leq k \leq n-1$, assume that players $j = 1, \dots, k$ eventually play good simple Smale plans.
 For a limit point $X = $ \\ $(x_1, \dots, x_k, z_{k+1}, \dots, z_n)$
 let $\bar x$ be the mean payoff among the good players and let $\bar z$ be the mean payoff among the dissenting players.
 If $\bar z \geq p_n$, then $X = (p_n, \dots, p_n)$ the complete cooperation point. Otherwise,
 $\bar x < \bar z < p_n$.  In particular, if $k = n-1$ then $z_n \geq p_n$ only when $X = (p_n, \dots, p_n)$.

 Now assume that $k \leq n-2$. If $z_{k+1} \geq p_n$ and $\tilde z$ is the mean payoff among the remaining dissenting players, then
 either $X = (p_n, \dots, p_n)$ or $\bar x > \tilde z$.
\end{theo}

{\bf Proof:} Because the strategies are assumed to be good simple Smale plans, it follows that $1 > \hat \l > \frac{n-1}{n}$.
Hence, $n \hat \l - k > 0$. It follows from (\ref{3.03}) that if $\bar z\geq p_n$, then $\bar x\geq p_n$. So $y\geq p_n$
%PZ $p_n - \bar x > 0$, then $p_n - \bar z > 0$. Otherwise,$\bar x = \bar z = p_n$ and so $y = p_n$
and this implies $X$ is the cooperation point. From the second equation in (\ref{3.03}) we
see that the dissenters on average do better than the good players but not as well as if they switched to cooperation and achieved the
$p_n$ payoff.

If there is a single dissenter, player $n$ then $\bar z = z_n$ and so $z_n \geq p_n$ only when $X = (p_n, \dots, p_n)$.

The second result follows from (\ref{3.06}) because $1 - n(1 - \hat \l)> 0$.

$\Box$ \vspace{.5cm}

Thus, if even one player uses a good simple Smale plan and the dissenting players do not allow the group to reach the cooperation point,
then on average they do worse than the $p_n$ payoff although they %the dissenters
do better than the good players.

With $\bar z < p_n$ it may still happen that some among the dissenters can reach a payoff of $p_n$ or greater. However if this happens then the
remaining dissenters on average do worse than the good players and so it pays for some of them to switch to a good simple Smale plan.
That is, at least one
of the remaining dissenters would do better by switching to a good  simple Smale plan.

The case of a dissenter obtaining a payoff greater than $p_n$ can happen. It seems to require exploiting
naive behavior by the other dissenters.
Here is an extreme case.

\begin{ex}\label{examplextra2}  With $n \geq 3$, if,  from some time on, player $k$ for $k = 1,\dots, n-2$ always cooperates, and player $n$
always defects, then player $n$ always receives a payoff between $r_{n-2}$ and $r_{n-1}$
regardless of the behavior of player $n-1$, e.g. even if player $n-1$ uses a good simple Smale strategy.
In particular, if $r_{n-2} > p_n$, then player $n$ always receives a payoff greater than $p_n$.

If, in addition, player $n-1$ eventually uses a simple Smale strategy then we obtain convergence to a unique
point $X$ regardless of the early plays.\end{ex}

{\bf Proof:} The strategies $All-C$ and $All-D$, always cooperating and always defecting, respectively,
trivially include Smale plans. If players $1, \dots, n-2$ always use $All-C$ and player $n$ always uses $All-D$, then
every outcome $S^T \in \O$ is either $V_0 = (p_{n-2}, \dots, p_{n-2}, r_{n-2}, r_{n-2})$ or $V_1 = (p_{n-1}, \dots, p_{n-1}, p_{n-1}, r_{n-1})$,
depending whether player $n-1$ defects or cooperates. It follows that every average point $s^T$ lies on the segment between them, which
we parameterize by $V_{a} = V_0 + a (V_1 - V_0)$ for $0 \leq a \leq 1$. If the strategies $All-C$ and $All-D$ are only adopted after some
time $T_0$, the effect of the initial terms on the average sequence $\{ s^T \}$ tends to $0$ as $T \to \infty$ and so in any case, the
limit set $\Omega$ is contained in the segment $[V_0,V_1]$. So for any $X \in \Omega$, $\pi_n(X) \geq r_{n-2}$.

The line $\d$ through $\Pi_{n-1}([V_0,V_1]) = [(r_{n-2},m_{n-2}),(p_{n-1},m_{n-1})]$ has negative slope. As a separation line $\ell$
has positive slope it follows that $\d \cap \ell$ is a single point.  Thus,  $\Pi_{n-1}^{-1}(\ell)$ is a hyperplane in $\R^n$
which intersects $[V_0,V_1]$ in a unique point $X$.  By Theorem \ref{theo3.2} $\Omega = \{ X \}$ and so we have convergence to $X$.

$\Box$ \vspace{.5cm}

 In the above example, suppose $\ell_{\l}$ is a separation line through $(p_n,p_n)$ with equation $y - p_n = \l (x - p_n)$ so that $\l \leq 1$.
We then let $X_{\l} = V_{a_{\l}}$ denote the
intersection point with the line $\d$. The value $a_{\l}$ is obtained by solving the equation $\bar \pi(V_{a}) - p_n = \l (\pi_{n-1}(V_{a}) - p_n)$.

In the special case $\l = 1$, $\ell_{1}$ is the diagonal with equation $y = x$ and so $a_1$ satisfies
\begin{equation}\label{3.07}
 m_{n-2} + a_1(m_{n-1} - m_{n-2}) = r_{n-2} +  a_1(p_{n-1} - r_{n-2}).
 \end{equation}
That is, $a_1$ equals
\begin{equation}\label{3.08}
\frac{r_{n-2} - m_{n-2}}{(m_{n-1} - m_{n-2}) + (r_{n-2} - p_{n-1})} = \frac{r_{n-2} - m_{n-2}}{(r_{n-2} - m_{n-2}) + (m_{n-1} - p_{n-1})}
\end{equation}

Thus, $\pi_n(X_1) = r_{n-2} + a_1(r_{n-1} - r_{n-2}) \geq p_n$ when
$a_1 \cdot r_{n-1} + (1 - a_1) \cdot r_{n-2} \geq p_n $.

If $\l < 1$, then $\ell_{\l}$ intersects $\d$ between $V_{a_1}$ and $V_1$. So if $\l < 1$, then $a_{\l} > a_1$. Because
$\pi_n(V_a)$ is increasing in $a$, we see that
\begin{equation}\label{3.09}
\l < 1 \quad \Longrightarrow \quad \pi_n(X_{\l}) > \pi_n(X_1).
\end{equation}

With $n = 3$, $(p_1, r_0, p_2, p_3, r_1, r_2) = (0, 2, 4, 6, 7, 8)$ provides an example with $r_1 > p_3$, since
$3 (m_0, m_1, m_2, m_3) = (6, 14, 16, 18)$. Thus, if player $1$ always cooperates and player $3$ always defects, then
player $3$ receives a payoff of at least $7 > 6$
regardless of the behavior of player $2$.

Even without the $r_{n - 2} > p_n$ assumption, it can happen that a dissenter obtains greater than $p_n$ even when one of
the other players uses a good, simple, Smale plan.

 \begin{ex}\label{examplextra3} Assume that for $n \geq 3$, $(p_1,r_0,p_2,r_1, \dots) = (0,1,2,3, \dots)$ as in Example \ref{exsimple}.
  If, eventually, player $k$ for $k = 1,\dots, n-2$ always cooperates, player $n-1$ uses a good, simple Smale plan and player $n$
always defects, then there is a unique limit point $X = (x_1, \dots, x_n)$  with
$x_n > p_n = 2n - 2$. \end{ex}

{\bf Proof:} In this case, we have
\begin{align}\label{3.10}
\begin{split}
\pi_j(V_a) = 2n - 6 + 2a = 2n - 3 \ + \ &(2a -3) \quad \text{ for } \quad j = 1, \dots, n-2, \\
\pi_{n-1}(V_a) = 2n - 3 - a, \qquad &\pi_{n}(V_a) = 2n - 3 + 2a, \\
\bar \pi(V_a) = 2n - 3 \  + \  &\frac{1}{n}[ a(2n - 3) - 3(n - 2)].
\end{split}
\end{align}
Equating $ \pi_{n-1}(V_a) = \bar \pi(V_a)$ to solve for $a_1$ we obtain $a_1 = \frac{n-2}{n-1}$. Hence,
\begin{equation}\label{3.11}
\pi_n(X_1) \ = \ 2n - 3 + 2 a_1 \ = \ 2n - 2 \ + \ \frac{n - 3}{n - 1}.
\end{equation}

Since $p_n = 2n - 2$, (\ref{3.11}) and (\ref{3.09}) imply $ \pi_n(X_{\l})  > \pi_n(X_1) \geq p_n$ with the latter inequality strict if
$n > 3$. Recall that for a good, simple Smale plan, the separation line passes through $(p_n,p_n)$ and has slope less than $1$.

$\Box$ \vspace{.5cm}

{\bfseries Remark:}  In the above example with $n = 3$ we easily see that if player $2$ uses the simple Smale plan with separation line the diagonal,
i.e. $\l = 1$,  then $a_1 = \frac{1}{2}$ and
the payoff vector is $X_1 = (1,2.5,4)$ when player $1$ uses $All-C$ and player $3$ uses $All-D$.  In general, if player $2$ uses the
diagonal simple Smale plan, then no player can obtain a payoff greater than $4 = p_3$, the cooperative payoff. For if $X = (x_1,x_2,x_3)$ is a limit point
then $\Pi_2(X) = (x_2, \bar x)$ lies on the diagonal by  Theorem \ref{theo3.2} and so $x_2 = \bar x$.
This implies that
%PZ We propose to replace the remaining arguments
%$\frac{1}{2}[(x_1,\bar x) + (x_3, \bar x)]  = (\bar x, \bar x)$.
%It follows that the furthest right one of these points
%can be is on the vertical line $x = 4$ since the furthest left the other can be is on the segment $[(0,2),(4,4)]$ which contains $(2,3)$.
%Observe that the equation of the line containing the latter segment is $x = 2y - 4$ and averaging we see $x = \frac{1}{2}[(2y - 4) + (4)] = y$.
%by the following
$\bar x=\frac{1}{2}(x_1+x_3)$. Observe that $\Pi_1(X)=(x_1,\,\bar x)\in\Pi_1(\O)= \S'$ which is the convex hull of
%$\{(0,2),\,(4,4),\,(1,1),\,(5,3)\}$. 
$\{(1,1),(0,2),(3,2),(2,3), (5,3), (4,4)\}$. It follows that
 $x_1\geq 2 \bar x-4$. Thus,  $x_3 = 2 \bar x - x_1 \leq 4$.

\vspace{1cm}

\bibliographystyle{amsplain}

\end{document}